# STATISTICAL ASPECTS OF THE FRACTIONAL STOCHASTIC CALCULUS

By Ciprian A. Tudor and Frederi G. Viens[1]

*Université de Paris 1 Panthéon–Sorbonne and Purdue University*

We apply the techniques of stochastic integration with respect to fractional Brownian motion and the theory of regularity and supremum estimation for stochastic processes to study the maximum likelihood estimator (MLE) for the drift parameter of stochastic processes satisfying stochastic equations driven by a fractional Brownian motion with any level of Hölder-regularity (any *Hurst* parameter). We prove existence and strong consistency of the MLE for linear and nonlinear equations. We also prove that a version of the MLE using only discrete observations is still a strongly consistent estimator.

**1. Introduction.** Stochastic calculus with respect to fractional Brownian motion (fBm) has recently experienced intensive development, motivated by the wide array of applications of this family of stochastic processes. For example, recent work and empirical studies have shown that traffic in modern packet-based high-speed networks frequently exhibits fractal behavior over a wide range of time scales; in quantitative finance and econometrics, the fractional Black–Scholes model has recently been introduced (see, e.g., [14, 17]) and this motivates the statistical study of stochastic differential equations governed by fBm.

The topic of parameter estimation for stochastic differential equations driven by standard Brownian motion is of course not new. Diffusion processes are widely used for modeling continuous time phenomena; therefore, statistical inference for diffusion processes has been an active research area over the last few decades. When the whole trajectory of the diffusion can be observed, then the parameter estimation problem is relatively simple, but of practical contemporary interest is work in which an approximate estimator,

Received March 2005; revised September 2006.

[1]Supported in part by NSF Grant DMS-02-04999.

using only information gleaned from the underlying process in discrete time, is able to do so as well as an estimator that uses continuously gathered information. Several methods have been employed to construct good estimators for this challenging question of discretely observed diffusions; among these methods, we refer to numerical approximation of the likelihood function (see [1, 5, 32]), martingale estimating functions (see [6]), indirect statistical inference (see [16]), the Bayesian approach (see [15]), some sharp probabilistic bounds on the convergence of estimators in [7], and [10, 12, 31] for particular situations. We mention the survey [36] for parameter estimation in discrete cases, further details in [21, 25] and the book [23].

Parameter estimation questions for stochastic differential equations driven by fBm are, in contrast, in their infancy. Some of the main contributions include [18, 19, 20, 33]. We take up these estimation questions in this article. Our purpose is to contribute further to the study of the statistical aspects of the fractional stochastic calculus by introducing the systematic use of efficient tools from stochastic analysis, to yield results which hold in some nonlinear generality. We consider the stochastic equation

$$(1) \qquad X_t = \theta \int_0^t b(X_s)\,ds + B_t^H, \qquad X_0 = 0,$$

where $B^H$ is a fBm with *Hurst* parameter $H \in (0,1)$ and the nonlinear function $b$ satisfies some regularity and nondegeneracy conditions. We estimate the parameter $\theta$ on the basis of the observation of the trajectory of the process $X$. The parameter $H$, which is assumed to be known, characterizes the local behavior of the process, with Hölder-regularity increasing with $H$; if $H = 1/2$, fBm is standard Brownian motion (BM) and thus has independent increments; if $H > 1/2$, the increments of fBm are positively correlated and the process is more regular than BM; if $H < 1/2$, the increments are negatively correlated and the process is less regular than BM. $H$ also characterizes the speed of decay of the correlation between distant increments. Estimating long-range dependence parameters is a difficult problem in itself which has received various levels of attention depending on the context; the text [3] can be consulted for an overview of the question; we have found the yet unpublished work [11], available online, which appears to propose a good solution applicable directly to fBm. Herein, we do not address the Hurst parameter estimation issue.

The results we prove in this paper are as follows: for *every $H$* in $(0,1)$,

- we give concrete assumptions on the nonlinear coefficient $b$ to ensure existence of the maximum likelihood estimator (MLE) for $\theta$ (Proposition 1);
- under certain hypotheses on $b$ which include nonlinear classes, we prove the strong consistency of the MLE (Theorems 2 and 3, depending on whether $H < 1/2$ or $H > 1/2$; and Proposition 2 and Lemma 3 for the scope of nonlinear applicability of these theorems); note that for $H > 1/2$ and $b$ linear, this has also been proved in [18];

- the bias and mean-square error for the MLE are estimated in the linear case (Proposition 3); this result was established for $H > 1/2$ in [18].

In this paper we also present a first practical implementation of the MLE studied herein, using only discrete observations of the solution $X$ of equation (1), by replacing integrals with their Riemann sum approximations. We show that

- in linear and some nonlinear cases, the discretization time-step for the Riemann sum approximations of the MLE can be fixed while still allowing for a strongly consistent estimator (Proposition 5 and Theorem 4).

To establish all these results, we use techniques in stochastic analysis including the *Malliavin* calculus, and supremum estimations for stochastic processes. The Malliavin calculus, or the stochastic calculus of variations, was introduced by P. Malliavin in [27] and developed by D. Nualart in [29]. Its original purpose was to study the existence and the regularity of the density of solutions to stochastic differential equations. Since our hypotheses in the present paper to ensure existence and strong consistency of the MLE are given in terms of certain densities [see condition (C)], the techniques of the Malliavin calculus appear as a natural tool.

We believe our paper is the first instance where the Malliavin calculus and supremum estimations are used to treat parameter estimation questions for fractional stochastic equations. These techniques should have applications and implications in statistics and probability reaching beyond the question of MLE for fBm. Indeed, even in the case of (Itô-) diffusion models, the strong consistency of an estimator follows if one can prove that an expression of the type $I_t := \int_0^t f^2(X_s)\,ds$ tends to $\infty$ as $t \to \infty$ almost surely, but a limited number of methods has been employed to deal with this kind of problem (if $X$ is Gaussian, the Laplace transform can be computed explicitly to show that $\lim_{t\to\infty} I_t = \infty$ a.s.; if $X$ is an ergodic diffusion, a local time argument can be used; particular situations have also been considered in [21, 22]). Our stochastic analytic tools constitute a new possibility, judging by the fact that the case of $H < 1/2$ is well within the reach of our tools, in contrast to the other above-mentioned methods, as employed, in particular, in [18] (see, however, a general Bayesian-type problem discussed in [19]).

The organization of our paper is as follows. Section 2 contains preliminaries on the fBm. In Section 3 we show the existence of the MLE for the parameter $\theta$ in (7) and in Section 4 we study its asymptotic behavior. Section 5 contains some additional results in the case when the drift function is linear. In Section 6 a discretized version of the MLE is studied. Some crucial technical proofs appear in the Appendix, while other more minor ones can be found in an extended version of this article posted on arXiv.org.

**2. Preliminaries on fBm and fractional calculus.** We consider $(B_t^H)_{t \in [0,T]}$, $B_0^H = 0$, a fBm with Hurst parameter $H \in (0,1)$, in a probability space $(\Omega, \mathcal{F}, \mathbf{P})$, that is, a centered Gaussian process with covariance function $R$ given by

$$(2) \qquad R(t,s) = \mathbf{E}[B_t^H B_s^H] = \tfrac{1}{2}(t^{2H} + s^{2H} - |t-s|^{2H}), \qquad s, t \in [0,T].$$

Let us denote by $K$ the kernel of the fBm such that (see, e.g., [28])

$$(3) \qquad B_t^H = \int_0^t K(t,s) \, dW_s,$$

where $W$ is a Wiener process (standard Brownian motion) under $\mathbf{P}$. Denote by $\mathcal{E}_H$ the set of step functions on $[0,T]$ and let $\mathcal{H}$ be the canonical Hilbert space of the fBm; that is, $\mathcal{H}$ is the closure of $\mathcal{E}$ with respect to the scalar product $\langle 1_{[0,t]}, 1_{[0,s]} \rangle_{\mathcal{H}} = R(t,s)$. The mapping $1_{[0,t]} \to B_t^H$ can be extended to an isometry between $\mathcal{H}$ and the Gaussian space generated by $B^H$ and we denote by $B^H(\varphi)$ the image of $\varphi \in \mathcal{H}$ by this isometry.

We also introduce the operator $K^*$ from $\mathcal{E}_H$ to $L^2([0,T])$ defined by

$$(4) \qquad (K^*\varphi)(s) = K(T,s)\varphi(s) + \int_s^T (\varphi(r) - \varphi(s)) \frac{\partial K}{\partial r}(r,s) \, dr.$$

With this notation, we have $(K^* 1_{[0,t]})(s) = K(t,s)$ and hence, the process

$$(5) \qquad W_t = \int_0^t (K^{*,-1} 1_{[0,t]})(s) \, dB_s^H$$

is a Wiener process (see [2]); in fact, it is *the* Wiener process referred to in formula (3), and for any nonrandom $\varphi \in \mathcal{H}$, we have $B^H(\varphi) = \int_0^T (K^*\varphi)(s) \, dW(s)$, where the latter is a standard Wiener integral with respect to $W$.

Last, we will use some elements of fractional calculus, with notation that can be found in [30] or our arXiv extended version, including the Riemann–Liouville fractional integral $I_{0+}^\alpha f(t)$ and derivative $D_{0+}^\alpha f(t)$ of any $f \in L^1[0,T]$ for $\alpha > 0$. The linear isomorphism $K_H$ from $L^2([0,T])$ onto $I_+^{H+1/2}(L^2([0,T]))$ whose kernel is $K(t,s)$ has an inverse which, depending on whether $H < 1/2$ or $H > 1/2$, equals

$$(6) \qquad \begin{aligned} (K_H^{-1}h)(s) &= s^{H-1/2} I_{0+}^{1/2-H}(s^{1/2-H} h'(s))(s), \text{ or} \\ &= s^{H-1/2} D_{0+}^{H-1/2}(s^{1/2-H} h'(s))(s). \end{aligned}$$

**3. The maximum likelihood estimator for fBm-driven stochastic differential equations.** We will analyze the estimation of the parameter $\theta \in \Theta \subset \mathbb{R}$ based on the observation of the solution $X$ of the stochastic differential equation

$$(7) \qquad X_t = \theta \int_0^t b(X_s)\,ds + B_t^H, \qquad X_0 = 0,$$

where $B^H$ is a fBm with $H \in (0, 1)$ and $b \colon \mathbb{R} \to \mathbb{R}$ is a measurable function. There are some strong known results concerning equation (7). In [30] strong existence and uniqueness is proved assuming only the linear growth

$$|b(x)| \le C(1 + |x|)$$

for $b$ when $H < 1/2$, and assuming Hölder-continuity of order $\alpha \in (1 - \frac{1}{2H}, 1)$ when $H > 1/2$. An extension is obtained in [8] when $H > 1/2$ if one adds a bounded nondecreasing left (or right) continuous function to the Hölder-continuous function $b$. In this paper we will only be using Lipschitz-continuous functions $b$.

Throughout the paper, we will typically avoid the use of explicit $H$-dependent constants appearing in the definitions of the operator kernels related to this calculus, since our main interest consists of asymptotic properties for estimators. In consequence, we will use the notation $C(H), c(H), c_H, \ldots$ for generic constants depending on $H$, which may change from line to line.

Our MLE construction is based on the following observation (see [30]). Consider the process $\tilde{B}_t^H = B_t^H + \int_0^t u_s\,ds$, where the process $u$ is adapted and with integrable paths. Then we can write

$$(8) \qquad \tilde{B}_t^H = \int_0^t K(t, s)\,dZ_s,$$

where

$$(9) \qquad Z_t = W_t + \int_0^t K_H^{-1}\left(\int_0^{\cdot} u_r\,dr\right)(s)\,ds.$$

We have the following Girsanov theorem.

THEOREM 1. (i) *Assume that $u$ is an adapted process with integrable paths such that*

$$t \to \int_0^t u_s\,ds \in I_{0+}^{H+1/2}(L^2([0, T])) \qquad a.s.$$

(ii) *Assume that $\mathbf{E}(V_T) = 1$, where*

$$(10) \qquad \begin{aligned} V_T = \exp\Big(&-\int_0^T K_H^{-1}\left(\int_0^{\cdot} u_r\,dr\right)(s)\,dW_s \\ &-\tfrac{1}{2}\int_0^T \left(K_H^{-1}\left(\int_0^{\cdot} u_r\,dr\right)(s)\right)^2 ds\Big). \end{aligned}$$

*Then under the probability measure $\tilde{\mathbf{P}}$ defined by $d\tilde{\mathbf{P}}/d\mathbf{P} = V_T$, it holds that the process $Z$ defined in* (9) *is a Brownian motion and the process $\tilde{B}^H$* (8) *is a fractional Brownian motion on $[0, T]$.*

HYPOTHESIS. We need to make, at this stage and throughout the remainder of the paper, the following assumption on the drift: $b$ is differentiable with bounded derivative $b'$; thus, the affine growth condition holds.

This Girsanov theorem is the basis for the following expression of the MLE.

PROPOSITION 1. *Denote, for every $t \in [0, T]$,*

$$(11) \qquad Q_t = Q_t(X) = K_H^{-1}\left(\int_0^{\cdot} b(X_r)\,dr\right)(t).$$

*Then $Q \in L^2([0, T])$ almost surely and the MLE is given by*

$$(12) \qquad \theta_t = -\frac{\int_0^t Q_s\,dW_s}{\int_0^t Q_s^2\,ds}.$$

Before proving Proposition 1, we need the following estimates.

LEMMA 1. *For every $s, t \in [0, T]$,*

$$(13) \qquad \sup_{s \leq t}|X_s| \leq \left(C\theta t + \sup_{s \leq t}|B_s^H|\right)e^{C\theta t}$$

*and*

$$(14) \qquad |X_t - X_s| \leq C\theta\left(1 + \sup_{u \leq T}|X_u|\right)|t - s| + |B_t^H - B_s^H|.$$

PROOF. With $C$ the linear growth constant of $b$, we have, for any $s$,

$$|X_s| \leq \int_0^s |\theta||b(X_u)|\,du + |B_s^H| \leq C\theta\int_0^s (1 + |X_u|)\,du + \sup_{u \leq s}|B_u^H|$$

and by Gronwall's lemma,

$$|X_s| \leq \left(C\theta s + \sup_{u \leq s}|B_u^H|\right)e^{C\theta s}, \qquad s \in [0, T],$$

and the estimate (13) follows. The second estimate follows by $b$'s affine growth. $\square$

PROOF OF PROPOSITION 1. Let

$$h(t) = \int_0^t b(X_s)\,ds.$$

We prove that the process $h$ satisfies (i) and (ii) of Theorem 1. Note first that the application of the operator $K_H^{-1}$ preserves the adaptability. We treat separately the cases when $H$ is greater than $1/2$ and less than $1/2$ .

*The case $H < 1/2$.* To prove (i), we only need to show that $Q \in L^2([0,T])$ **P**-a.s. Now using relation (6) we thus have, for some constant $C_H$ which may change from line to line, using the hypothesis $|b(x)| \leq C(1 + |x|)$, for all $s \leq T$,

$$
(15) \quad
\begin{aligned}
|Q_s| &\leq C_H s^{H-1/2} \left| \int_0^s (s-u)^{-1/2-H} u^{1/2-H} b(X_u) \, du \right| \\
&\leq C_H \left( 1 + \sup_{u \leq s} |X_u| \right),
\end{aligned}
$$

which we can rewrite, thanks to Lemma 1, as

$$
\sup_{s \leq T} |Q_s| \leq C(H,T) \left( 1 + \sup_{s \leq T} |X(s)| \right),
$$

which, thanks to inequality (13), is of course much stronger than $Q \in L^2([0,T])$ a.s.

To prove (ii), it suffices to show that there exists a constant $\alpha > 0$ such that

$$
\sup_{s \leq T} \mathbf{E}(\exp(\alpha Q_s^2)) < \infty.
$$

Since $Q$ satisfies (15), the above exponential moment is a trivial consequence of inequality (13) and Fernique's theorem on the exponential integrability of the square of a seminorm of a Gaussian process.

*The case $H > 1/2$.* Using formula (6), we have in this case that

$$
(16) \quad
\begin{aligned}
Q_s &= c_H \Bigg[ s^{1/2-H} b(X_s) \\
&\quad + \left( H - \frac{1}{2} \right) s^{H-1/2} \int_0^s \frac{b(X_s) s^{1/2-H} - b(X_u) u^{1/2-H}}{(s-u)^{H+1/2}} \, du \Bigg] \\
&= c_H \Bigg[ s^{1/2-H} b(X_s) \\
&\quad + \left( H - \frac{1}{2} \right) s^{H-1/2} b(X_s) \int_0^s \frac{s^{1/2-H} - u^{1/2-H}}{(s-u)^{H+1/2}} \, du \\
&\quad + \left( H - \frac{1}{2} \right) s^{H-1/2} \int_0^s \frac{b(X_s) - b(X_u)}{(s-u)^{H+1/2}} u^{1/2-H} \, du \Bigg]
\end{aligned}
$$

and using the fact that

$$\int_0^s (s^{1/2-H} - u^{1/2-H})(s-u)^{-H-1/2}\,du = c(H)s^{1-2H},$$

we get

$$|Q_s| \leq c_H \left( s^{1/2-H}|b(X_s)| + s^{H-1/2}\int_0^s \frac{b(X_s) - b(X_u)}{(s-u)^{H+1/2}} u^{1/2-H}\,du \right)$$
$$:= A(s) + B(s).$$

The first term $A(s)$ above can be treated as in [30], proof of Theorem 3, due to our Lipschitz assumption on $b$. We obtain that, for every $\lambda > 1$,

$$(17) \qquad \mathbf{E}\left( \exp\left( \lambda \int_0^t A_s^2\,ds \right) \right) < \infty.$$

To obtain the same conclusion for the second summand $B(s)$, we note that, by Lemma 1, up to a multiplicative constant, the random variable

$$G = \sup_{0 \leq u < s \leq T} \frac{|X_s - X_u|}{|u - s|^{H-\varepsilon}}$$

is bounded by

$$\left( 1 + \sup_{u \leq T}|X_u| \right)|t - s|^{1-H+\varepsilon} + \sup_{0 \leq u < s \leq T} \frac{|B_s^H - B_u^H|}{|u - s|^{H-\varepsilon}}$$

and it suffices to use the calculations contained in [30].

*Conclusion.* Properties (i) and (ii) are established for both cases of $H$, and we may apply Theorem 1. Expression (12) for the MLE follows a standard calculation, since [using the notation $P_\theta$ for the probability measure induced by $(X_s)_{0 \leq s \leq t}$ and the fact that $P_0 = \mathbf{P}$]

$$(18) \qquad F(\theta) := \log \frac{dP_\theta}{dP_0} = -\theta \int_0^t Q_s\,dW_s - \frac{\theta^2}{2}\int_0^t Q_s^2\,ds. \qquad \square$$

We finish this section with some remarks that will relate our construction to previous work [18, 19, 33]. Details about these links are given in Section 5.

*Alternative form of the MLE.* By (7), we can write, by integrating the quantity $K^{*,-1}1_{[0,t]}(s)$ for $s$ between 0 and $t$,

$$(19) \qquad \int_0^t (K^{*,-1}1_{[0,t]}(\cdot))(s)\,dX_s = \theta \int_0^t (K^{*,-1}1_{[0,t]}(\cdot))(s)b(X_s)\,ds + W_t.$$

On the other hand, by (7) again,

$$(20) \qquad X_t = \int_0^t K(t,s)\,dZ_s,$$

where $Z$ is given by (9). Therefore, we have the equality

$$(21) \qquad \int_0^t (K^{*,-1}1_{[0,t]}(\cdot))(s)\,dX_s = Z_t.$$

By combining (19) and (21), we obtain

$$\int_0^t K_H^{-1}\left(\int_0^{\cdot} b(X_r)\,dr\right)(s)\,ds = \int_0^t (K^{*,-1}1_{[0,t]}(\cdot))(s)b(X_s)\,ds$$

and thus, the function $t \to \int_0^t (K^{*,-1}1_{[0,t]}(\cdot))(s)b(X_s)\,ds$ is absolutely continuous with respect to the Lebesgue measure and

$$(22) \qquad Q_t = \frac{d}{dt}\int_0^t (K^{*,-1}1_{[0,t]}(\cdot))(s)b(X_s)\,ds.$$

By (9), we get that the function (18) can be written as

$$F(\theta) = -\theta\int_0^t Q_s\,dZ_s + \frac{\theta^2}{2}\int_0^t Q_s^2\,ds.$$

As a consequence, the maximum likelihood estimator $\theta_t$ has the equivalent form

$$(23) \qquad \theta_t = \frac{\int_0^t Q_s\,dZ_s}{\int_0^t Q_s^2\,ds}.$$

*The above formula* (23) *shows explicitly that the estimator $\theta_t$ is observable if we observe the whole trajectory of the solution $X$.*

**4. Asymptotic behavior of the maximum likelihood estimator.** This section is devoted to studying the strong consistency of the MLE (12). A similar result has been proven in the case $b(x) \equiv x$ and $H > 1/2$ in [18]. We propose here a proof of strong consistency for a class of functions $b$ which contains significant nonlinear examples. By replacing (9) in (23), we obtain that

$$\theta_t - \theta = \frac{\int_0^t Q_s\,dW_s}{\int_0^t Q_s^2\,ds},$$

with $Q$ given by (11) or (22). To prove that $\theta_t \to \theta$ almost surely as $t \to \infty$ (which means by definition that the estimator $\theta_t$ is strongly consistent), by the strong law of large numbers we need only show that

$$(24) \qquad \lim_{t\to\infty}\int_0^t Q_s^2\,ds = \infty \qquad \text{a.s.}$$

To prove that $\lim_{t\to\infty} \int_0^t Q_s^2\, ds = \infty$ in a nonlinear case, it is necessary to make some assumption of nondegeneracy on the behavior of $b$. In order to illustrate our method using the least amount of technicality, we will restrict our study to the case where the function $|b|$ satisfies a simple probabilistic estimate with respect to fractional Brownian motion:

(C) There exist positive constants $t_0$ and $K_b$, both depending only on $H$ and the function $b$, and a constant $\gamma < 1/(1+H)$ such that, for all $t \geq t_0$ and all $\varepsilon > 0$, we have $\tilde{\mathbf{P}}[|Q_t(\tilde{\omega})|/\sqrt{t} < \varepsilon] \leq \varepsilon t^{\gamma H} K_b$, where under $\tilde{\mathbf{P}}$, $\tilde{\omega}$ has the law of fractional Brownian motion with parameter $H$.

4.1. *The case $H < 1/2$.* In this part we prove the following result.

THEOREM 2. *Assume that $H < 1/2$ and that condition* (C) *holds. Then the estimator $\theta_t$ is strongly consistent, that is,*

$$\lim_{t\to\infty} \theta_t = \theta \qquad almost\ surely.$$

Before proving this theorem, we discuss condition (C). To understand this condition, we first note that, with $\mu_H^t$ the positive measure on $[0,t]$ defined by $\mu_H^t(dr) = (r/t)^{1/2-H}(t-r)^{-1/2-H}\, dr$, according to the representation (6) we have

$$Q_t = \int_0^t \mu_H^t(ds) b(\tilde{\omega}_s)$$

and therefore, by the change of variables $r = s/t$,

$$(25) \qquad \frac{Q_t}{\sqrt{t}} = \int_0^1 \mu_H^1(dr)\frac{b(\tilde{\omega}_{tr})}{t^H}$$

$$(26) \qquad \overset{\mathcal{D}}{=} \int_0^1 \mu_H^1(dr)\frac{b(t^H\tilde{\omega}_r)}{t^H},$$

where the last inequality is in distribution under $\tilde{\mathbf{P}}$.

If $b$ has somewhat of a linear behavior, we can easily imagine that $b(t^H\tilde{\omega}_r)/t^H$ will be of the same order as $b(\tilde{\omega}_r)$. Therefore, $Q_t/\sqrt{t}$ should behave, in distribution for fixed $t$, similarly to the universal random variable $\int_0^1 \mu_H^1(dr)b(\tilde{\omega}_r)$ (whose distribution depends only on $b$ and $H$). Generally speaking, if this random variable has a bounded density, the strongest version of condition (C), that is, with $\gamma = 0$, will follow. In the linear case, of course, the factors $t^H$ disappear from expression (26), leaving a random variable which is indeed known to have a bounded density, uniformly in $t$, by the arcsine law. The presence of the factor $t^{\gamma H}$ in condition (C) gives even more flexibility, however, since, in particular, it allows a bound on the density of $Q_t/\sqrt{t}$ to be proportional to $t^H$.

Leaving aside these vague considerations, we now give, in Proposition 2, a simple sufficient condition on $b$ which implies condition (C). The proof of this condition uses the tools of the Malliavin calculus; as such, it requires some extra regularity on $b$. We also give a class of nonlinear examples of $b$'s satisfying (C) [condition (28) in Lemma 3] which are more restricted in their global behavior than in Proposition 2, but do not require any sort of local regularity for $b$.

PROPOSITION 2. *Assume $H < 1/2$. Assume that $b'$ is bounded and that $b''$ satisfies $|b''(x)| \leq b_1/(1 + |x|^\beta)$ for some $\beta \in (H/(1-H), 1)$. Assume that $|b'|$ is bounded below by a positive constant $b_0$. Then, letting $\gamma = 1 - \beta$, condition* (C) *holds.*

REMARK 1. The condition $\gamma < 1/(1 + H)$ from condition (C) translates as $\beta > H/(1-H)$, which is consistent with $\beta < 1$ because $H < 1/2$.

The nondegeneracy condition on $|b'|$ above can be relaxed. If, for $x \geq x_0$, $|b'(x)| \geq x^{-\alpha}$ holds, then condition (C) holds as long as $\alpha$ does not exceed a positive constant $\alpha_0(H)$ depending only on $H$. We omit the very technical proof of this fact.

The hypothesis of fractional power decay on $b''$, while crucial, does allow $b$ to have a highly nonlinear behavior. Compare with Lemma 3 below, which would correspond to the case $\beta = 1$ here.

The hypotheses of the above proposition imply that $b$ is monotone.

The proof of Proposition 2 requires a criterion from the Malliavin calculus, which we present here. The book [29] by D. Nualart is an excellent source for proofs of the results we quote. Here we will only need to use the following properties of the Malliavin derivative $D$ with respect to $W$ [recall that $W$ is the standard Brownian motion used in the representation (3), i.e., defined in (5)]. For simplicity of notation, we assume that all times are bounded by $T = 1$. The operator $D$, from a subset of $L^2(\Omega)$ into $L^2(\Omega \times [0, 1])$, is essentially the only one which is consistent with the following two rules:

1. Consider a centered random variable in the Gaussian space generated by $W$ (first chaos); it can therefore be represented as $Z = W(f) = \int_0^1 f(s) \, dW(s)$ for some nonrandom function $f \in L^2([0, 1])$. The operator $D$ picks out the function $f$, in the sense that, for any $r \in [0, 1]$,

$$D_r Z = f(r).$$

2. $D$ is compatible with the chain rule, in the sense that, for any $\Phi \in C^1(\mathbf{R})$ such that both $F := \Phi(Z)$ and $\Phi'(Z)$ belongs to $L^2(\Omega)$, for any $r \in [0, 1]$,

$$D_r F = D_r \Phi(Z) = \Phi'(Z) D_r Z = \Phi'(Z) f(r).$$

For instance, using these two rules, definition (3) and formula (5) relative to the fBm $\tilde{\omega}$ under $\tilde{\mathbf{P}}$, we have that, under $\tilde{\mathbf{P}}$, for any $r \leq s$,

$$(27) \qquad D_r b(t^H \tilde{\omega}_s) = t^H b'(t^H \tilde{\omega}_s) K(s, r).$$

It is convenient to define the domain of $D$ as the subset $\mathbf{D}^{1,2}$ of r.v.'s $F \in L^2(\Omega)$ such that $D.F \in L^2(\Omega \times [0,1])$. Denote the norm in $L^2([0,1])$ by $\|\cdot\|$. The set $\mathbf{D}^{1,2}$ forms a Hilbert space under the norm defined by

$$\|F\|_{1,2}^2 = \mathbf{E}|F|^2 + \mathbf{E}\|D.F\|^2 = \mathbf{E}|F|^2 + \mathbf{E}\int_0^1 |D_r F|^2 \, dr.$$

Similarly, we can define the second Malliavin derivative $D^2 F$ as a member of $L^2(\Omega \times [0,1]^2)$, using an iteration of two Malliavin derivatives, and its associated Hilbert space $\mathbf{D}^{2,2}$. Non-Hilbert spaces, using powers other than 2, can also be defined. For instance, the space $\mathbf{D}^{2,4}$ is that of random variables $F$ having two Malliavin derivatives and satisfying

$$\|F\|_{2,4}^4 = \mathbf{E}|F|^4 + \mathbf{E}\|D.F\|^2 + \mathbf{E}\|D^2_{\cdot,\cdot}F\|_{L^2([0,1]^2)}^4$$

$$= \mathbf{E}|F|^4 + \mathbf{E}\int_0^1 |D_r F|^2 \, dr + \mathbf{E}\left(\int_0^1 |D_r D_s F|^2 \, dr \, ds\right)^2 < \infty.$$

We also note that the so-called *Ornstein–Uhlenbeck* operator $L$ acts as follows (see [29], Proposition 1.4.4):

$$LF = L\Phi(Z) = -Z\Phi'(Z) + \Phi''(Z)\|f\|^2.$$

We have the following Lemma, whose proof we omit because it follows from ([29], Proposition 2.1.1 and Exercise 2.1.1).

LEMMA 2. *Let $F$ be a random variable in $\mathbf{D}^{2,4}$ such that $\mathbf{E}[\|DF\|^{-8}] < \infty$. Then $F$ has a continuous and bounded density $f$ given by*

$$f(x) = \mathbf{E}\left[1_{(F>x)}\left(\frac{-LF}{\|DF\|^2} - 2\frac{\langle DF \otimes DF; D^2 F\rangle_{L^2([0,1]^2)}}{\|DF\|^4}\right)\right].$$

PROOF OF PROPOSITION 2. Using (26), and the notation $\mu = \mu_1^H$, let

$$F = \frac{Q_t}{\sqrt{t}} = \int_0^1 \mu(dr) \frac{b(t^H \tilde{\omega}_r)}{t^H}.$$

It is sufficient to prove that $F$ has a density which is bounded by $K_b t^{\gamma H}$, where the constant $K_b$ depends only on $b$ and $H$. Indeed, $\tilde{\mathbf{P}}[|Q_t(\tilde{\omega})|/\sqrt{t} < \varepsilon] \leq \int_0^\varepsilon K_b t^{\gamma H} \, dx = \varepsilon t^{\gamma H} K_b$. In this proof, $C_{b,H}$ denotes a constant depending only on $b$ and $H$, whose value may change from line to line.

*Step* 1: *calculating the terms in Lemma* 2. We begin with the calculation of $DF$. Since the Malliavin derivative is linear, we get $D_r F = t^{-H} \int_0^1 \mu(ds) \times D_r(b(t^H \tilde{\omega}_s))$. Then from (27) we get $D_r F = \int_r^1 \mu(ds) b'(t^H \tilde{\omega}_s) K(s, r)$. Thus, we can calculate

$$\|DF\|^2 = \int_0^1 dr \left| \int_r^1 \mu(ds) b'(t^H \tilde{\omega}_s) K(s, r) \right|^2$$

$$= \int_0^1 \int_0^1 \mu(ds) \mu(ds') b'(t^H \tilde{\omega}_s) b'(t^H \tilde{\omega}_{s'}) \int_0^{\min(s, s')} K(s, r) K(s, r') \, dr$$

$$= \int_0^1 \int_0^1 \mu(ds) \mu(ds') b'(t^H \tilde{\omega}_s) b'(t^H \tilde{\omega}_{s'}) R(s, s'),$$

where $R$ is the covariance of fBm in (2). A similar calculation yields

$$D_{q,r}^2 F = t^H \int_{\max(q,r)}^1 \mu(ds) b''(t^H \tilde{\omega}_s) K(s, r) K(s, q)$$

and

$$\|D^2 F\|_{L^2([0,1]^2)}^2 = t^{2H} \int_0^1 \int_0^1 \mu(ds) \mu(ds') b''(t^H \tilde{\omega}_s) b''(t^H \tilde{\omega}_{s'}) |R(s, s')|^2.$$

For the Ornstein–Uhlenbeck operator, which is also linear, we get

$$-LF = \int_0^1 \mu(ds)(b'(t^H \tilde{\omega}_s) \tilde{\omega}_s + b''(t^H \tilde{\omega}_s) t^H s^{2H}).$$

*Step* 2: *estimating the terms in Lemma* 2. With the expressions in the previous step, using the hypotheses of the proposition, we now obtain for some constant $C_H$ depending only on $H$,

$$\tilde{\mathbf{E}}[|LF|] \leq C_H \left( \|b'\|_\infty + t^H b_1 \int_0^1 \mu(ds) s^{2H} \tilde{\mathbf{E}} \left[ \frac{1}{1 + t^{\beta H} |\tilde{\omega}_s|^\beta} \right] \right)$$

$$= C_H \left( \|b'\|_\infty + t^H b_1 \mathbf{E} \left[ \int_0^1 \mu(ds) s^{2H} \frac{1}{1 + (ts)^{\beta H} |Z|^\beta} \right] \right),$$

where $Z$ is a generic standard normal random variable. We deal first with the integral in $s$. For $s \in [0, 1/2]$, $\mu(ds)$ has a bounded density, and thus, for any $a > 0$ and any $\alpha < 1$, we have $\int_0^{1/2} ds (1 + as^\alpha)^{-1} \leq (1 - \alpha)^{-1} a^{-1}$; on the other hand, for $s \in [1/2, 1]$, we can bound $(1 + (ts)^{\beta H} |Z|^\beta)^{-1}$ above by $(t/2)^{-\beta H} |Z|^{-\beta-1}$. We immediately obtain, using $a = t^{\beta H} |Z|^\beta$ and $\alpha = \beta H$,

$$\tilde{\mathbf{E}}[|LF|] \leq C_{H,b} \left( 1 + t^{H(1-\beta)} \frac{1}{1 - \beta H} \mathbf{E}[|Z|^{-\beta}] \right.$$

$$\left. + t^{H(1-\beta)} \mathbf{E}[|Z|^{-\beta}] \int_{1/2}^1 \frac{ds}{(1-s)^{1/2+H}} \right)$$

$$\leq C_{H,b} (1 + t^{H(1-\beta)}).$$

The estimation of $\|D^2 F\|$ is similar. Using its expression in the previous step, the measure $d\mu/ds \cdot d\mu/ds' \cdot R(s, s')$ and the fact that $\int_0^1 ds(1 + as^\alpha)^{-2} \leq (1 - 2\alpha)^{-1}a^{-2}$, with $\alpha = \beta H < 1/2$, we get

$$\tilde{\mathbf{E}}\|D^2 F\|_{L^2([0,1]^2)} \leq C_{b,H} t^{H(1-\beta)}.$$

Also, almost surely, for any $p \geq 2$, for some constant $C_{H,p}$ depending only on $H$ and $p$, since $b'$ has a constant sign, we obtain

$$\frac{1}{\|DF\|^p} = \left( \int_{[0,1]^2} \mu(ds)\mu(ds')R(s, s')|b'(t^H\tilde{\omega}_s)b'(t^H\tilde{\omega}_{s'})| \right)^{-p/2} \leq C_{H,p}b_0^{-p}.$$

Last, it is convenient to invoke the Cauchy–Schwarz inequality to get

$$\frac{\langle DF \otimes DF; D^2 F \rangle_{L^2([0,1]^2)}}{\|DF\|^4} \leq \frac{\|D^2 F\|_{L^2([0,1]^2)}}{\|DF\|^2}.$$

*Step* 3: *applying Lemma* 2; *conclusion.*  The third estimate in the previous step (for $p = 8$) proves trivially that $\tilde{\mathbf{E}}\|DF\|^{-8}$ is finite. That $F \in \mathbf{D}^{2,4}$ follows again trivially from the boundedness of $b'$ and $b''$ using the expressions in Step 1. Thus, Lemma 2 applies. We conclude from the estimates in the previous step that $F$ has a density $f$ which is bounded as

$$f(x) \leq C_{H,b}(1 + t^{H(1-\beta)})b_0^{-2}.$$

With $t \geq 1$, the conclusion of the proposition follows.  □

A smaller class of functions $b$ satisfying condition (C) and covering all $H \in (0, 1)$ is given in the following, proved in the extended version of this article on arXiv.org.

LEMMA 3.  *Let* $H \in (0, 1)$. *Assume* $xb(x)$ *has a constant sign for all* $x \in \mathbf{R}_+$ *and a constant sign for all* $x \in \mathbf{R}_-$. *Assume*

$$(28) \qquad\qquad |b(x)/x| = c + h(x)$$

*for all* $x$, *where* $c$ *is a fixed positive constant and* $\lim_{x \to \infty} h(x) = 0$. *Then condition* (C) *is satisfied with* $\gamma = 0$.

Condition (C) also holds for any $b$ of the above form to which a constant $C$ is added: $|(b(x) - C)/x| = c + h(x)$ and $\lim_{x \to \infty} h(x) = 0$. Note that this condition is less restrictive than saying $b$ is asymptotically affine, since it covers the family $b(x) = C + cx + (|x| \wedge 1)^\alpha$ for any $\alpha \in (0, 1)$. In some sense, condition (C) with $\gamma = 0$ appears to be morally equivalent to this class of functions.

PROOF OF THEOREM 2.

*Step* 1: *setup.* Since we only want to show that (24) holds, and since $\int_0^t |Q_s|^2\,ds$ is increasing, it is sufficient to satisfy condition (24) for $t$ tending to infinity along a sequence $(t_n)_{n\in\mathbf{N}}$. We write, according to (6), for each fixed $t \geq 0$,

$$I_t = I_t(X) := \int_0^t |Q_s(X)|^2\,ds = \int_0^t \left| \int_0^s \mu_H^s(dr) b(X_r) \right|^2\,ds,$$

where $X$ is the solution of the Langevin equation (7) and the positive measure $\mu_H^s$ is defined by $\mu_H^s(dr) = (r/s)^{1/2-H}(s-r)^{-1/2-H}\,dr$. Since the theorem's conclusion is an almost sure statement about $X$, and from the Girsanov Theorem 1 applied to $X$ the measures $\mathbf{P}$ and $\tilde{\mathbf{P}}$ are equivalent, it is sufficient to work under $\tilde{\mathbf{P}}$, that is to assume that $X = \tilde{\omega}$ is a standard fBm; we omit the dependence of $Q_s$ and $I_t$ on $X = \tilde{\omega}$ for simplicity.

We will show that with $t_n = n^k$ for some positive integer $k$ chosen below, $I_{t_n}$ converges to $\infty$, by restricting the integration defining $I_{t_n}$ to a small interval of length $b_n = n^{-j}$ near $t_n$, where $j$ will be another integer chosen below. Indeed,

$$
\begin{aligned}
I_{t_n} = \int_0^{t_n} |Q_s|^2\,ds &\geq \int_{t_n-b_n}^{t_n} |Q_s|^2\,ds \\
&\geq b_n |Q_{t_n}|^2 - \int_{t_n-b_n}^{t_n} |Q_{t_n} - Q_s||Q_{t_n} + Q_s|\,ds \\
&\geq b_n \left( |Q_{t_n}|^2 - \sup_{s\in[t_n-b_n,t_n]} |Q_{t_n} - Q_s||Q_{t_n} + Q_s| \right) \\
&\geq b_n |Q_{t_n}|^2 - 2b_n \sup_{s\in[0,t_n]} |Q_s| \sup_{s\in[t_n-b_n,t_n]} |Q_{t_n} - Q_s| \\
&:= A_n - B_n.
\end{aligned}
$$

(29)

*Step* 2: *the diverging term $A_n$.* The term $A_n$ is easily shown to converge to infinity almost surely thanks to condition (C) modulo a condition on $j$ and $k$. Indeed, we have, for any sequence $a_n$,

$$\mathbf{P}[b_n |Q_{t_n}|^2 < a_n^2] \leq K_b t_n^{\gamma H} a_n (b_n t_n)^{-1/2} = K_b n^{k(\gamma H-1/2)+j/2} a_n.$$

In order for the right-hand side of the last expression to be summable in $n$ while being able to choose $\lim_{n\to\infty} a_n = 0$, it is sufficient to impose

$$(30) \qquad\qquad j + 2 < (1 - 2\gamma H)k.$$

Thus, specifically, under this condition, with $a_n^2 = n^\ell$ for $0 < \ell < (1-2\gamma H)k - j - 2$, we have by the Borel–Cantelli lemma the existence of an almost surely finite $n_0$ such that $n > n_0$ implies $A_n = b_n |Q_{t_n}|^2 > n^\ell$.

*Step* 3: *the error term* $B_n$. To control the term $B_n$, we need two estimates, whose proofs are based on elementary facts on the trajectories of fBm, namely, boundedness and Hölder-continuity, which are derived from that process's Gaussian and scaling properties. Details are in this article's extended arXiv.org version. For any sequences $b_n, t_n$ such that $b_n \ll t_n$ (meaning $\lim_n b_n/t_n = 0$), for any fixed $\varepsilon > 0$, there exists an almost-surely finite random variable $c_\varepsilon$ such that, for all $n > 0$, we have

$$(31) \qquad \sup_{s \in [0, t_n]} |Q_s| \leq c_\varepsilon (t_n)^{1/2 + \varepsilon}$$

and

$$(32) \qquad \sup_{s, t \in [t_n - b_n, t_n]} |Q_t - Q_s| \leq c_\varepsilon (t_n)^{1/2 - H + \varepsilon} (b_n)^{H - \varepsilon}.$$

One immediately obtains that almost-surely

$$B_n \leq 2 c_\varepsilon^2 b_n (t_n)^{1/2 + \varepsilon} (t_n)^{1/2 - H + \varepsilon} (b_n)^{H - \varepsilon} = 2 c_\varepsilon^2 n^{k(1 - H + 2\varepsilon) - j(1 + H - \varepsilon)}.$$

The statement that $B_n$ converges to 0 (i.e., that $\varepsilon$ can be chosen positive while having the power in the last expression above be negative) now follows from assuming that $j$ and $k$ are related by

$$(33) \qquad k(1 - H) - j(1 + H) < 0.$$

*Step* 4: *conclusion.* We may now use the results of the last two parts, namely, that $A_n > n^\ell$ while $\lim_{n \to \infty} B_n = 0$, to conclude from (29) that the statement of the theorem holds, provided that conditions (30) and (33) hold, that is,

$$\frac{j + 2}{1 - 2\gamma H} < k < \frac{1 + H}{1 - H} j.$$

The theorem now follows because the relation $\gamma < 1/(1 + H)$ in condition (C) implies that if $j$ is large enough, we do indeed have $\frac{j+2}{1 - 2\gamma H} < \frac{1+H}{1-H} j$. □

4.2. *The case* $H > 1/2$. Due to the fact that the function $Q$ is less regular in this case, we should not expect the proof of the following theorem to be a consequence of the proof in the case $H < 1/2$. Nevertheless, it deviates from the former proof very little. On the other hand, we cannot rely on Proposition 2 to find a convenient sufficient condition for condition (C); instead, we can look to the nonlinear class of examples in Lemma 3, which satisfy the strong version ($\gamma = 0$) of condition (C) for all $H \in (0, 1)$. The next result's proof is in the Appendix.

THEOREM 3. *Assume that $H > 1/2$ and $b$ satisfies condition* (C) *with $\gamma = 0$ [e.g., $b$ satisfies condition* (28) *in Lemma* 3]. *Then the maximum likelihood estimator $\theta_t$ is strongly consistent.*

**5. The linear case.** In this section we present some comments in the case when the drift $b$ is linear. We will assume that $b(x) \equiv x$ to simplify the presentation. In this case, the solution $X$ to equation (7) is the fractional Ornstein–Uhlenbeck process and it is possible to prove more precise results concerning the asymptotic behavior of the maximum likelihood estimator.

REMARK 2. In [9] it is shown that, for any $H \in (0, 1)$, there is a unique almost surely continuous process $X$ satisfying (7), and it can be represented as

$$(34) \qquad X_t = \int_0^t e^{\theta(t-u)} \, dB_u^H, \qquad t \in [0, T],$$

where the above integral is a Wiener integral with respect to $B^H$ (which exists also as a pathwise Riemann–Stieltjes integral). It follows from the stationarity of the increments of $B^H$ that $X$ is stationary and the decay of its autocovariance function is like a power function. The process $X$ is ergodic, and for $H > 1/2$ it exhibits long-range dependence.

Let us briefly recall the method employed in [18] to estimate the drift parameter of the fractional OU process. Let us consider the function, for $0 < s < t \leq 1$,

$$k(t, s) = c_H^{-1} s^{1/2 - H} (t - s)^{1/2 - H} \qquad \text{with } c_H = 2H\Gamma(\tfrac{3}{2} - H)\Gamma(H + \tfrac{1}{2})$$

(35)

and let us denote its Wiener integral with respect to $B^H$ by

$$(36) \qquad M_t^H = \int_0^t k(t, s) \, dB_s^H.$$

It has been proved in [28] that $M^H$ is a Gaussian martingale with bracket

$$(37) \qquad \langle M^H \rangle_t := \omega_t^H = \lambda_H^{-1} t^{2-2H} \qquad \text{with } \lambda_H = \frac{2H\Gamma(3 - 2H)\Gamma(H + 1/2)}{\Gamma(3/2 - H)}.$$

The authors called $M^H$ the fundamental martingale associated to fBm. Therefore, observing the process $X$ given by (7) is the same thing as observing the process

$$Z_t^{KB} = \int_0^t k(t, s) \, dX_s,$$

which is actually a semimartingale with the decomposition

$$(38) \qquad Z_t^{KB} = \theta \int_0^t Q_s^{KB} \, d\omega_s^H + M_t^H,$$

where

$$(39) \qquad Q_t^{KB} = \frac{d}{d\omega^H} \int_0^t k(t, s) X_s \, ds, \qquad t \in [0, T].$$

By using Girsanov's theorem (see [18] and [28]), we obtain that the MLE is given by

$$\theta_t := \theta_t^{KB} = \frac{\int_0^t Q_s^{KB} \, dZ_s^{KB}}{\int_0^t (Q_s^{KB})^2 \, d\omega_s^H}. \tag{40}$$

REMARK 3.   We can observe that our operator (12) or (23) coincides (possibly up to a multiplicative constant) with the one used in [18] and given by (40). Assume that $H < 1/2$; the case $H > 1/2$ is just a little more technical.

PROOF.   Using relations (11) and (35), we can write

$$Q_t = C(H)t^{H-1/2} \int_0^t s^{1/2-H}(t-s)^{-1/2-H} b(X_s) \, ds$$

$$= C(H)t^{H-1/2} \int_0^t \frac{d}{dt} k(t,s) b(X_s) \, ds$$

$$= C(H)t^{H-1/2} \frac{d}{dt} \int_0^t k(t,s) b(X_s) \, ds.$$

It is not difficult to see that $\frac{d}{dt} \int_0^t k(t,s) b(X_s) \, ds = C(H)t^{1-2H} Q_t^{KB}$ and therefore,

$$Q_t = C(H)t^{1/2-H} Q_t^{KB}. \tag{41}$$

On the other hand, it can be similarly seen that

$$Z_t^{KB} = C(H) \int_0^t s^{1/2-H} \, dZ_s \tag{42}$$

and the estimations given by (40) and (23) coincide up to a constant.   □

To compute the expression of the bias and of the mean square error and to prove the strong consistency of the estimator, one has the option, in this explicit linear situation, to compute the Laplace transform of the quantity $\int_0^t (Q_s^{KB})^2 \, d\omega_s^H$. This is done for $H > 1/2$ in [18], Section 3.2, and the following properties are obtained:

- the estimator $\theta_t$ is strongly consistent, that is,

$$\theta_t \to \theta \qquad \text{almost surely when } t \to \infty;$$

- the *bias* and the *mean square error* are given by:
  - If $\theta < 0$, when $t \to \infty$, then

$$\mathbf{E}(\theta_t - \theta) \backsim \frac{2}{t}, \qquad \mathbf{E}(\theta_t - \theta)^2 \backsim \frac{2}{t}|\theta|; \tag{43}$$

– If $\theta > 0$, when $t \to \infty$, then

$$(44) \qquad \mathbf{E}(\theta_t - \theta) \curlywedge -2\sqrt{\pi \sin \pi H}\theta^{3/2}e^{-\theta t}\sqrt{t},$$

$$(45) \qquad \mathbf{E}(\theta_t - \theta)^2 \curlywedge 2\sqrt{\pi \sin \pi H}\theta^{5/2}e^{-\theta t}\sqrt{t}.$$

It is interesting to realize that the rate of convergence of the bias and of the mean square error does not depend on $H$. In fact, the only difference between the classical case (see [26]) and the fractional case is the presence of the constant $\sqrt{\pi H}$ in (43), (44) and (45). It is natural to expect the same results if $H < 1/2$. This is true, as stated next, and proved in Section A.2.

PROPOSITION 3. *If $H < 1/2$, then* (43), (44) *and* (45) *hold.*

**6. Discretization.** In this last section we present a discretization result which allows the implementation of an MLE for an fBm-driven stochastic differential equation.

We first provide background information on discretely observed diffusion processes in the classical situation when the driving noise is the standard Brownian motion. Assume that

$$dX_t = b(X_t, \theta) + \sigma(X_t, \theta)\,dW_t,$$

where $\sigma, b$ are known functions, $W$ is a standard Wiener process and $\theta$ is the unknown parameter. If continuous information is available, the parameter estimation by using the maximum likelihood method is somewhat simpler; indeed, the maximum likelihood function can be obtained by means of the standard Girsanov theorem and there are results on the asymptotic behavior of the estimator (consistency, efficiency etc...). We refer to the monographs [4, 34] or [24] for complete expositions of this topic.

"Real-world" data is, however, typically discretely sampled (e.g., stock prices collected once a day or, at best, at every tick). Therefore, statistical inference for discretely observed diffusions is of great interest for practical purposes and at the same time it poses a challenging problem. Here the main obstacle is the fact that discrete-time transition functions are not known analytically and consequently, the likelihood function is in general not tractable. In this situation there are alternative methods to treat the problem. Among these methods, we refer to numerical approximation to the likelihood function (see [1, 5, 32]), martingale estimating functions (see [6]), indirect statistical inference (see [16]) and Bayesian approaches (see [15]). We refer to [36] for a survey of methods of estimation in the discrete case. When the transition functions of the diffusion $X$ are known, and $\sigma(x, \theta) = \sigma x$ with $\sigma$ unknown and not depending on $\theta$, then Dacunha–Castelle and Florens–Zmirou [12] propose a maximum likelihood estimator which is strongly consistent for the pair $(\theta, \sigma)$. They also give a measure of the loss of information due to the discretization as a function depending on the interval between two observations.

A more particular situation is the case when $\sigma$ is known (assume that $\sigma = 1$). Then the maximum likelihood function, given by $\exp(\theta \int_0^t b(X_s) \, dX_s - \frac{\theta^2}{2} \int_0^t b(X_s)^2 \, ds)$, can been approximated using Riemann sums as

$$\exp\left( \theta \sum_{i=0}^{N-1} b(X_{t_i})(X_{t_{i+1}} - X_{t_i}) - \frac{\theta^2}{2} \sum_{i=0}^{N-1} b(X_{t_i})^2 (t_{i+1} - t_i) \right).$$

As a consequence, the following maximum likelihood estimator can be obtained from the discrete observations of the process $X$ at times $t_0, \ldots, t_N$ in a fixed interval $[0, T]$, with discrete mesh size decreasing to 0 as $N \to \infty$:

$$(46) \qquad \theta_{N,T} = -\frac{\sum_{i=0}^{N-1} b(X_{t_i})(X_{t_{i+1}} - X_{t_i})}{\sum_{i=0}^{N-1} |b(X_{t_i})|^2 (t_{i+1} - t_i)}$$

(see [25], including proof of convergence to the continuous MLE).

In the fractional case, we are aware of no such results. We propose a first concrete solution to the problem. We choose to work with the formula (23) by replacing the stochastic integral in the numerator and the Riemann integral in the denominator by their corresponding approximate Riemann sums, using discrete integer time. Specifically, we define, for any integer $n \geq 1$,

$$(47) \qquad \bar{\theta}_n := \frac{\sum_{m=0}^{n} Q_m (Z_{m+1} - Z_m)}{\sum_{m=0}^{n} |Q_m|^2}.$$

Our goal in this section is to prove that $\bar{\theta}_n$ is in fact a consistent estimator for $\theta$. By our Theorems 2 and 3, it is of course sufficient to prove that $\lim_{n \to \infty} (\bar{\theta}_n - \theta_n) = 0$ almost surely. One could also consider the question of the discretization of $\theta_T$ using a fine time mesh for fixed $T$, and show that this discretization converges almost surely to $\theta_t$; by time-scaling, such a goal is actually equivalent to our own.

It is crucial to note that in the fractional case the process $Q$ given by (11) depends continuously on $X$ and, therefore, the discrete observation of $X$ does not allow one to obtain directly the discrete observation of $Q$. We explain how to remedy this issue: $Q_m$ appearing in (47) can be easily approximated if we know the values of $X_n, n \geq 1$, since only a deterministic integral appears in the expression of (11); indeed, for $H < 1/2$, the quantity

$$(48) \qquad \check{Q}_n = c(H) n^{H-1/2} \sum_{j=0}^{n-1} (n-j)^{-H-1/2} j^{1/2-H} b(X_j)$$

can be deduced from observations and it holds that $\lim_n (Q_n - \check{Q}_n) = 0$ almost surely. This last fact requires proof, which is simpler than the proof of convergence of $\bar{\theta}_n - \theta_n$ to 0, but still warrants care; we present the crucial estimates of this proof in Section A.3.1.

Note, moreover, that calculation of $\bar{\theta}_n$ also relies on $Z_m$, which is not observable; yet from formula (21), where $Z_m$ is expressed as a stochastic integral of a deterministic function against the increments of $X$, again, we may replace all the $Z_m$'s by their Riemann sum; proving that these sums converge to the $Z_m$'s follows from calculations which are easier than those presented in Section A.3.1, because they only require discretizing the deterministic integrand. We summarize this discussion in the following precise statement, referring to Section A.3.1 for indications of its proof.

PROPOSITION 4.    With $\check{Q}_n$ as in (48) and $\check{Z}_n = \sum_{j=0}^{n-1} (K^{*,-1} 1_{[0,n]}(\cdot))(j) \times (X_{j+1} - X_j)$, then almost surely $\bar{\theta}_n - \check{\theta}_n$ converges to 0, where $\check{\theta}_n$ is given by (47) with $Z$ and $Q$ replaced by $\check{Z}$ and $\check{Q}$.

Let $\langle M \rangle_n$ denote the quadratic variation at time $n$ of a square-integrable martingale $M$. We introduce the two semimartingales

$$(49) \qquad A_t := \int_0^t Q_s \, dZ_s,$$

$$(50) \qquad B_t := \int_0^t Q_{[s]} \, dZ_s,$$

where $[s]$ denotes the integer part of $s$. We clearly have $B_n = \sum_{m=0}^{n-1} Q_m (Z_{m+1} - Z_m)$. Thus, using the fact that $Z$ is a Brownian motion under $\bar{\mathbf{P}}$, we see that

$$(51) \qquad \langle B \rangle_n = \sum_{m=0}^{n-1} |Q_m|^2,$$

while

$$(52) \qquad \langle A - B \rangle_n = \int_0^n |Q_s - Q_{[s]}|^2 \, ds = \sum_{m=0}^{n-1} \int_m^{m+1} |Q_s - Q_m|^2 \, ds.$$

Therefore, from definitions (12) and (47), we immediately get the expressions

$$\theta_n = \frac{A_n}{\langle A \rangle_n} \quad \text{and} \quad \bar{\theta}_n = \frac{B_n}{\langle B \rangle_n}.$$

The following proposition defines a strategy for proving that $\bar{\theta}_n$ – and, by the previous proposition, $\check{\theta}_n$ – is a consistent estimator for $\theta$. See Section A.3.2 for its proof.

PROPOSITION 5.    Let $H \in (0,1)$. If there exists a constant $\alpha > 0$ such that

- $n^\alpha \langle A - B \rangle_n / \langle B \rangle_n$ is bounded almost surely for $n$ large enough,
- for all $k \geq 1$, for some constant $K > 0$, almost surely, for large $n$, $\langle B \rangle_n^k \geq K \mathbf{E}[\langle B \rangle_n^k]$,
- and for all $k > 1$, $\mathbf{E}[|\langle A - B \rangle_n|^k] \leq n^{-k\alpha} \mathbf{E}[|\langle B \rangle_n|^k]$,

then almost surely $\lim_{n \to \infty} \bar{\theta}_n = \theta$.

The following theorem, proved in Section A.3.3 under the condition (C′) below, which is stronger than (C), still allows for nonlinear examples.

THEOREM 4.   *Assume $b'$ is bounded and the following condition holds:*

(C′)   *There exist constant $t_0, K_b > 0$ depending only on $H$ and $b$, such that, for all $t \geq t_0$ and all $\varepsilon > 0$, $\tilde{\mathbf{P}}[|Q_t(\tilde{\omega})|/\sqrt{t} < \varepsilon] \leq \varepsilon K_b$, where under $\tilde{\mathbf{P}}$, $\tilde{\omega}$ has the law of fBm with parameter $H$.*

*Then for all $H \in (0, 1/2)$, almost surely $\lim_{n \to \infty} \bar{\theta}_n = \theta$, where the discretization $\bar{\theta}_n$ of the maximum likelihood estimator $\theta_n$ is defined in (47). If $H \in (1/2, 1)$, the same conclusion holds if we assume in addition that $b''$ is bounded.*

*By Proposition 4, the above statements hold with $\bar{\theta}$ replaced by $\check{\theta}$.*

REMARK 4.   Condition (C′) holds as soon as the random variable $Q_t(\tilde{\omega})/\sqrt{t}$ has a density that is bounded uniformly $t$. When $H < 1/2$, this is a statement about the random variables $\int_0^1 \mu_1^H(ds) b(t^H \tilde{\omega}_s) t^{-H}$. In all cases, condition (C′) holds for the class of nonlinear functions defined in Lemma 3.

We conjecture that Theorem 4 holds if we replace (C′) by (C), in view, for example, of the fact that the conditions of Proposition 5 hold for any $\alpha < 2H$. Step 1 in the theorem's proof is the obstacle to us establishing this.

## APPENDIX

**A.1. Proof of Theorem 3.**   Recall from the proof of Proposition 1 that we can write

$$(53) \qquad Q_t = c(H) t^{1/2-H} b(X_t) + c'(H) \int_0^t \mu_H^t(dr) b(X_t) - b(X_r)).$$

We note that in this case the expression $\mu_H^t(dr)$ does not determine a measure, but we still use this notation to simplify the presentation; the Lipschitz assumption on $b$ and the Hölder property of $X$ do ensure the existence of the integral.

One can actually follow the proof in the case $H < 1/2$ line-by-line. All we have to do here is prove an equivalent of relations (31) and (32) on the supremum and variations of $Q$, this being the only point where the form of $Q$, which differs depending on whether $H$ is greater or less than $1/2$, is used. We briefly indicate how the second summand of $Q$ in (53) (which is the most difficult to handle) can be treated.

Under $\tilde{\mathbf{P}}$, we denote $X$ by $\tilde{\omega}$, since its law is that of standard fBm. The quantity $Q'_t := \int_0^t \mu^t_H(dr)(b(X_t) - b(X_r))$ equals $\int_0^1 \mu^1_H(dr) t^{-H}(b(t^H \tilde{\omega}_1) - b(t^H \tilde{\omega}_r))$ in distribution. Now, let $V'_t := t^{-1/2} Q'_t$. Omitting the details, we state that instead of relations (31) and (32) on $Q'$, it is equivalent to show the following bound for some $M > 2$:

$$(54) \qquad \tilde{\mathbf{E}}\Big[\sup_{s,t \in [t_n - b_n, t_n]} |V'_t - V'_s|^M\Big] \le C_{M,H,b}\Big(\frac{b_n}{t_n}\Big)^{HM},$$

which follows from some elementary calculations, and a standard application of the Kolmogorov lemma on continuity (see [35], Theorem I.2.1). Details are in the extended version of this article on arXiv. A direct proof, via (31) and (32), not via (54), is also possible.

**A.2. Proof of Proposition 3.** To avoid tedious calculations with fractional integrals and derivatives, we will take advantage of the calculations performed in [18] when $H > 1/2$; nevertheless, we believe that a direct proof is also possible. Actually, the only step where the authors of [18] use the fact that $H$ is greater than $1/2$ is the computation of the process $Q$. By relations (20) and (22), we can write

$$Q_t = \frac{d}{dt} \int_0^t \int_v^t (K^{*,-1} 1_{[0,t]}(\cdot))(s) K(s,v)\, ds\, dZ_v.$$

Note that from the formulas presented in Section 2, we have

$$(K^{*,-1} 1_{[0,t]}(\cdot))(s) = c(H) s^{1/2-H} \int_s^t u^{1/2-H}(u-s)^{-H-1/2}, \qquad H < \tfrac{1}{2},$$

$$(K^{*,-1} 1_{[0,t]}(\cdot))(s) = c(H) s^{1/2-H} \frac{d}{ds} \int_s^t u^{1/2-H}(u-s)^{-H+1/2}, \qquad H > \tfrac{1}{2}.$$

To unify the notation, we write

$$(K^{*,-1} 1_{[0,t]}(\cdot))(s) = c(H) s^{1/2-H} \frac{d}{ds} \int_s^t u^{1/2-H}(u-s)^{-H+1/2}, \qquad H \in (0,1),$$

and we just observe that the constant $c(H)$ above is analytic with respect to $H$. Let us consider, for $v \le t$, a function $A(v,t)$ such that

$$\int_v^t A(v,s)\, ds = \int_v^t (K^{*,-1} 1_{[0,t]}(\cdot))(s) K(s,v)\, ds.$$

Then, obviously, $Q_t = \int_0^t A(t,v)\, dZ_v$.

On the other hand, it has been proved in [18] [see relations (3.4) and (3.5) therein] that, for $H > 1/2$, $Q_t^{KB} = \int_0^t A^{KB}(t,v)\, dZ_v^{KB}$ with $A^{KB}(t,s) =$ $c(H)(t^{2H-1} + s^{2H-1})$. Using the relations between $Q$ and $Q^{KB}$ and between $Z$ and $Z^{KB}$ (see Remark 3), it follows that, for every $H > 1/2$, and $s < t$,

$$(55) \qquad A(s,t) = c(H) \left[ \left( \frac{s}{t} \right)^{1/2-H} + \left( \frac{t}{s} \right)^{1/2-H} \right].$$

We show that the above relation (55) is true for $H < 1/2$ as well. We use an argument inspired by [13], proof of Theorem 3.1. We observe that the functions

$$H \in (0,1) \to A(s,t) \quad \text{and} \quad H \in (0,1) \to c(H) \left[ \left( \frac{s}{t} \right)^{1/2-H} + \left( \frac{t}{s} \right)^{1/2-H} \right]$$

are analytic with respect to $H$ and coincide on $(1/2, 1)$. Moreover, both are well-defined for every $H \in (0,1)$ (in fact, it follows from [18] that $A$ is well defined for $H > 1/2$ and it is more regular for $H \le 1/2$). To conclude (55) for every $H \in (0,1)$, we invoke the fact that if $f, g : (a,b) \to \mathbb{R}$ are two analytic functions and the set $\{x \in (a,b); f(x) = g(x)\}$ has an accumulation point in $(a,b)$, then $f = g$.

As a consequence, (55) holds for every $H \in (0,1)$ and this shows that

$$\int_0^t Q_s \, dZ_s = \int_0^t Q_s^{KB} \, dZ_s^{KB} = c(H) \left( Z_t^{KB} \int_0^t r^{2H-1} \, dZ_r^{KB} - t \right)$$

and all calculations contained in [18], Sections 3.2, 4 and 5 hold for every $H$. $\square$

### A.3. Proof of Theorem 4.

A.3.1. *Proof of Proposition 4.* For conciseness, we only indicate how to establish one of the crucial estimates for this proposition, that the quantity

$$S_n := \frac{\sum_{m=0}^n (Q_m - \check{Q}_m)(Z_{m+1} - Z_m)}{\sum_{m=0}^n |Q_m|^2}$$

converges to 0 almost surely, and then only for $H < 1/2$. Since we want to show that $S_n$ tends to 0 almost surely, and $\mathbf{P}$ and $\tilde{\mathbf{P}}$ share the same null sets, we may assume that $Z$ is a Brownian motion, and $X$ is a fractional Brownian motion adapted to $Z$'s filtration.

Define the quantity

$$R_n = \sum_{j=0}^{n-1} (n-j)^{-H-1/2} j^{1/2-H} \int_j^{j+1} (b(X_j) - b(X_s)) \, ds.$$

This is related to $S_n$ via the fact that $m^{H-1/2}|R_m| = Q_m - \check{Q}_m$. We claim that, for any $\varepsilon > 0$, almost surely, for large $m$, that is, $m \ge m_0$, $|R_m| \le$ $r_0 + m^{-H+1+\varepsilon} c_H \|b'\|$, where $r_0$ is a fixed random variable. This is sufficient to conclude that $\lim_n S_n = 0$. Indeed, we will see below [Section A.3.3, Step 2, inequality (61)] that, almost surely, for large $n$, $\sum_{m=0}^n |Q_m|^2 \geq n^2$. Then the sum of all terms in the numerator of $S_n$ for $m \leq m_0$, after having been divided by $S_n$'s denominator, tends to 0 when $n \to \infty$. On the other hand, the IID terms $\{Z_{m+1} - Z_m\}_{m \in \mathbf{N}}$ are standard normal, so that one trivially proves that almost surely for $n \geq m_0$ (abusively using the same $m_0$ as above), up to some nonrandom universal constant $c$, $|Z_{m+1} - Z_m| \leq c\sqrt{\log m}$. It follows that the portion of $S_n$ for $m \geq m_0$ is bounded above by $n^{-2} \sum_{m=0}^n m^{H-1/2}(r_0 + m^{-H+1+\varepsilon} c_H \|b'\|)\sqrt{\log m}$, which is itself bounded above by $(r_0 + c_H \|b'\|)n^{3/2+\varepsilon}$, which obviously tends to 0 as $n \to \infty$ as soon as $\varepsilon < 1/2$.

Now let us prove our claim on $R_m$. It is a known fact, which is obtained using standard tools from Gaussian analysis, or simply the Kolmogorov lemma (see [35], Theorem I.2.1), that, for any $M \geq 1$,

$$\mathbf{E}\left[\sup_{s,t \in [j,j+1]} |X_t - X_s|^M\right] \leq j^{HM}.$$

The usual application of the Borel–Cantelli lemma after Chebyshev's inequality for an $M$ large enough implies that, for any $\alpha > H$, almost surely, for large $j$, $\sup_{s,t \in [j,j+1]} |X_t - X_s| \leq j^\alpha$. Consequently, for any $\varepsilon > 0$,

$$\begin{aligned}
|R_m| &\leq 2\|b'\| \sum_{j=0}^{m_0} (n-j)^{-H-1/2} j^{1/2-H} \int_j^{j+1} |X_j - X_s|\, ds \\
&\quad + \|b'\| \sum_{j=m_0}^{n-1} (n-j)^{-H-1/2} j^{1/2-H} \sup_{s \in [j,j+1]} |X_j - X_s| \\
&= r_0 + n^{-2H}\|b'\| \sum_{j=m_0}^{n-1} (1-j/n)^{-H-1/2} (j/n)^{1/2-H} \sup_{s \in [j,j+1]} |X_j - X_s| \\
&\leq r_0 + n^{-2H} n^{H+\varepsilon} \sum_{j=m_0}^{n-1} (1-j/n)^{-H-1/2} (j/n)^{1/2-H} (j/n)^{H+\varepsilon} \\
&= r_0 + c_H n^{-H+\varepsilon+1}(1 + O(1/n)),
\end{aligned}$$

where the last estimate is by virtue of the Riemann sums for $\int_0^1 \frac{x^{1/2+\varepsilon}}{(1-x)^{H+1/2}}\, dx$.

A.3.2. *Proof of Proposition* 5. By our Theorems 2 and 3, it is of course sufficient to prove that $\lim_{n\to\infty}(\bar{\theta}_n - \theta_n) = 0$. In preparation for this, we first note that by classical properties for quadratic variations, and using our hypothesis, for large enough $n$ we have

$$
|\langle B\rangle_n - \langle A\rangle_n| = |\langle (B-A), (B+A)\rangle_n|
$$

(56)
$$
\leq |\langle B+A\rangle_n|^{1/2} |\langle B-A\rangle_n|^{1/2}
$$

$$
\leq \sqrt{2} n^{-\alpha} |\langle B\rangle_n|^{1/2} |\langle A\rangle_n + \langle B\rangle_n|^{1/2}.
$$

Now we prove that (56) implies, almost surely,

(57)
$$
\lim_{n\to\infty} \frac{\langle A\rangle_n}{\langle B\rangle_n} = 1.
$$

Indeed, let $x_n = \langle A\rangle_n / \langle B\rangle_n$. Then we can write

$$
|x_n - 1| = \frac{|\langle B\rangle_n - \langle A\rangle_n|}{\langle B\rangle_n}
$$

$$
\leq \sqrt{2} n^{-\alpha} |\langle B\rangle_n|^{-1/2} |\langle A\rangle_n + \langle B\rangle_n|^{1/2}
$$

$$
= c\sqrt{2} n^{-\alpha} |1 + x_n|^{1/2},
$$

where $c$ is a possibly random almost surely finite constant. Let $\varepsilon > 0$ be given; it is elementary to check that the inequality $(x-1)^2 \leq 2\varepsilon(x+1)$ is equivalent to $|x - (1+\varepsilon)| \leq \sqrt{4\varepsilon + \varepsilon^2}$. For us this implies immediately $|x_n - 1| \leq 6 c n^{-\alpha}$, proving the claim (57).

Now we have

(58)
$$
\theta_n - \bar{\theta}_n = \frac{A_n}{\langle A\rangle_n} - \frac{B_n}{\langle B\rangle_n} = \frac{A_n - B_n}{\langle B\rangle_n} + A_n \frac{\langle B\rangle_n - \langle A\rangle_n}{\langle A\rangle_n \langle B\rangle_n}.
$$

Using (56), we have that the second term in (58) is bounded above in absolute value by

$$
\sqrt{2} n^{-\alpha} \frac{A_n}{\langle A\rangle_n} \frac{|\langle A\rangle_n + \langle B\rangle_n|^{1/2}}{|\langle B\rangle_n|^{1/2}} = \sqrt{2} n^{-\alpha} \frac{A_n}{\langle A\rangle_n} \left( \frac{\langle A\rangle_n}{\langle B\rangle_n} + 1 \right)^{1/2}.
$$

By Theorems 2 and 3, $A_n / \langle A\rangle_n$ converges to the finite constant $\theta$. By the limit (57), the last term in the above expression converges to 2, so that the entire expression converges to 0. Let $k$ and $\gamma$ be fixed positive values. For the first term in (58), using our hypotheses, by the Chebyshev and Burkholder–Davis–Gundy inequalities, and from the expression of the semimartingales $Z$ as $Z_t = \int_0^t Q_s \, dW_s + \theta \int_0^t Q_s \, ds$, we have

$$
\mathbf{P}[|A_n - B_n|^k > n^{-k\gamma} \mathbf{E}[\langle B\rangle_n^k]] \leq n^{\gamma k} \mathbf{E}^{-1}[\langle B\rangle_n^k] \mathbf{E}[|A_n - B_n|^k]
$$

$$
\leq c(\theta) 2^k n^{\gamma k} n^{-k\alpha}.
$$

Thus, picking a positive value $\gamma < \alpha$ and choosing $k$ large enough, by the Borel–Cantelli lemma, almost surely, for $n$ large enough,

$$
|A_n - B_n| \leq n^{-\gamma} \mathbf{E}[\langle B\rangle_n^k]^{1/k} \leq \frac{1}{K} n^{-\gamma} \langle B\rangle_n,
$$

which finishes the proof of the proposition.

A.3.3. *Proof of Theorem* 4 (*Steps* 1 *through* 4). In this entire proof $n_0(\omega)$ will denote a random almost surely finite integer; it may change from line to line, as it is introduced via various different applications of the Borel–Cantelli lemma, but one only needs to take the supremum of all such integers to have correct statements throughout.

*Step* 0. *Strategy.* First note that since the probability measures $\mathbf{P}$ and $\tilde{\mathbf{P}}$ are equivalent (see Theorem 1), almost sure statements under one measure are equivalent to statements about the same stochastic processes under the other measure, and therefore, we may prove the statements in the theorem by assuming that the process $Z$ in the definitions (49) and (50) is a standard Brownian motion, since such is its law under $\tilde{\mathbf{P}}$. Furthermore, for the same reason, we can assume that, in these same definitions, $Q$ is given by formula (11), where $X$ is replaced by $\tilde{\omega}$ whose law is that of standard fBm. We will use specifically, instead of (11), the explicit formula (25) when $H < 1/2$. For $H > 1/2$, the formula (16) must be used instead, which shows the need for a control of $b$'s second derivative. For the sake of conciseness, we restrict our proof to the case $H < 1/2$. The result of the theorem is established as soon as one can verify the hypotheses of Proposition 5. Here we present only the proof of the first of the three hypotheses. The other two are proved using similar or simpler techniques. To achieve our goal in this proof, it is thus sufficient to prove that, almost surely, for large $n$, $\langle B \rangle_n \geq n^{\alpha_1}$, while $\langle A - B \rangle_n \leq n^{\alpha_2}$, where the values $\alpha_1$ and $\alpha_2$ are nonrandom and $\alpha_1 > \alpha_2$. We establish these estimates in the Appendix.

*Step* 1. *Bounding* $|Q|^2$ *below.* Using only condition (C'), we immediately get, for any $\gamma \in (0, 1/2 - H)$, for any large $t$,

$$\mathbf{P}[|Q_t| < t^{1/2 - \gamma}] \leq K_b t^{-\gamma}.$$

To be able to apply the Borel–Cantelli lemma, we now let $t = n^A$, where $n$ is an integer and $A$ is a constant exceeding $\gamma^{-1}$. We then get, almost surely, for any $n > n_0(\omega)$,

$$|Q_{n^A}| > n^{A(1/2 - \gamma)}. \tag{59}$$

We also bound other $Q_{m'}$'s that are in close proximity to $Q_{n^A}$. For any fixed integer $\bar{m}$, consider the set $I_{\bar{m}}$ of integers $m'$ in the interval $[\bar{m} - \bar{m}^{4\gamma}, \bar{m}]$, where $\gamma$ is also assumed to be less than $1/4$. Then by inequality (32), for any $\varepsilon > 0$, for some almost-surely finite r.v. $c_\varepsilon$,

$$\sup_{m' \in I_{\bar{m}}} |Q_{\bar{m}} - Q_{m'}| \leq c_\varepsilon \bar{m}^{1/2 - H + \varepsilon} (\bar{m}^{4\gamma})^{H - \varepsilon}$$

$$= c_\varepsilon \bar{m}^{1/2 - (H - \varepsilon)(1 - 4\gamma)}$$

$$= c_\varepsilon \bar{m}^{1/2 - \varepsilon'},$$

where $\varepsilon' = (H - \varepsilon)(1 - 4\gamma)$ is positive for $\varepsilon < H$ and our hypothesis on $\gamma$. Thus with $0 < \varepsilon' < H(1 - 4\gamma)$, almost surely, we may write that

$$|Q_{m'}| > |Q_{\bar{m}}| - c_\varepsilon \bar{m}^{1/2 - \varepsilon'}.$$

Certainly, if $\bar{m}$ is of the form $n^A$ for large enough $n$, by choosing $\gamma$ small enough, we obtain that the lower bound $\bar{m}^{1/2 - \gamma}$ on $|Q_{\bar{m}}|^2$ obtained in (59) is dominant compared to $\bar{m}^{1/2 - \varepsilon}$ for $\varepsilon$ close to $H(1 - 4\gamma)$. Hence, we get, almost surely, with $\bar{m} = n^A$ large enough, for all $m' \in [\bar{m} - \bar{m}^{4\gamma}, \bar{m}]$,

$$(60) \qquad |Q_{m'}|^2 > |Q_{\bar{m}}|^2 \left( 1 - \frac{c_\varepsilon \bar{m}^{1/2 - \varepsilon'}}{|Q_{\bar{m}}|} \right)^2 \geq |Q_{\bar{m}}|^2 / 2.$$

*Step* 2. *Bounding* $\langle B \rangle$ *from below.* For $n$ given, let $n_1$ be the largest integer such that $n_1^A \leq n < (n_1 + 1)^A$. Also assume $n$ is large enough so that $n_1^A \geq n_0(\omega)$. Thus, applying (60) with $\bar{m} = n_1^A$,

$$\langle B \rangle_m \geq \sum_{m=0}^{n_1^A} |Q_m|^2 \geq |Q_{n_1^A}|^2 + \sum_{m' = n_1^A - (n_1^A)^{4\gamma}} |Q_{m'}|^2$$

$$\geq |Q_{n_1^A}|^2 + 2^{-1} \sum_{m' = n_1^A - (n_1^A)^{4\gamma}} |Q_{n_1^A}|^2 \geq |Q_{n_1^A}|^2 (n_1^A)^{4\gamma}.$$

We can now invoke (59) to say that, almost surely, for $n > n_0(\omega)^A$,

$$\langle B \rangle_n \geq (n_1)^{A(2 - 2\gamma)} (n_1^A)^{4\gamma} = 2^{-1} (n_1)^{A(2 + 2\gamma)}.$$

Given that we may write $n_1^A(1 + n_1^{-1}) > n$, so that $n_1^A > n/2$, we can finally conclude that

$$(61) \qquad \langle B \rangle_n \geq \frac{1}{2^{1 + 2A(1 + \gamma)}} n^{2 + 2\gamma}.$$

*Step* 3. *Bounding* $\langle A - B \rangle$*'s terms from above.* We may generically bound the general term of $\langle A - B \rangle_n$ directly using the bound (32) and its associated random variable $c_\varepsilon$: almost surely, for all $m \geq 0$,

$$\int_m^{m+1} |Q_s - Q_m|^2 \, ds \leq \int_m^{m+1} ds \sup_{t \in [m, m+1]} |Q_t - Q_m|^2$$

$$(62) \qquad = \sup_{t \in [m, m+1]} |Q_t - Q_m|^2$$

$$\leq c_\varepsilon (m^{1/2 - H + \varepsilon} 1^{H - \varepsilon})^2 = m^{1 - 2H + 2\varepsilon}.$$

We conclude that, for any $\delta = 1 - \varepsilon < H$, almost surely, for all $m$,

$$\int_m^{m+1} |Q_s - Q_m|^2 \, ds \leq c_\varepsilon m^{1 - 2\delta}.$$

*Step* 4. *Conclusion.* From the formula $\langle A - B \rangle_n = \sum_{m=0}^{n-1} \int_m^{m+1} |Q_s - Q_m|^2 \, ds$, using the last estimate of the previous step, we get

$$\langle A - B \rangle_n \leq c_\varepsilon n^{2-2\delta}.$$

From the final estimate (61) of Step 2, we may now write almost surely

$$\frac{\langle A - B \rangle_n}{\langle B \rangle_n} \leq \frac{c_\varepsilon}{n^{2(\gamma+\delta)}}.$$

Hence, the first statement of Proposition 5 is established for any $\alpha < 2H$.

**Acknowledgments.** We gratefully acknowledge our debt to the insightful comments of the Editor, Associate Editor and two referees, which resulted in several important improvements on an earlier version of this paper.

Samos–Matisse
Centre d'Economie de La Sorbonne
Université de Paris 1 Panthéon–Sorbonne
90 rue de Tolbiac
75634 Paris
France
E-mail: tudor@univ-paris1.fr

Department of Statistics
  and Department of Mathematics
Purdue University
150 N. University St.
West Lafayette, Indiana 47907-2067
USA
E-mail: viens@purdue.edu